%% Plain TeX
%%
%% Numerotation automatatique sous plain: Ch. Sorger
\input amssym.def
\newcount\secno
\newcount\prmno
\newif\ifnotfound
\newif\iffound
\def\namedef#1{\expandafter\def\csname #1\endcsname}
\def\nameuse#1{\csname #1\endcsname}

\long\def\ifundefined#1#2#3{\expandafter\ifx\csname
  #1\endcsname\relax#2\else#3\fi}
\def\hwrite#1#2{{\let\the=0\edef\next{\write#1{#2}}\next}}

% Working with lists
\toksdef\ta=0 \toksdef\tb=2
\long\def\leftappenditem#1\to#2{\ta={\\{#1}}\tb=\expandafter{#2}%
                                \edef#2{\the\ta\the\tb}}
\long\def\rightappenditem#1\to#2{\ta={\\{#1}}\tb=\expandafter{#2}%
                                \edef#2{\the\tb\the\ta}}

\def\lop#1\to#2{\expandafter\lopoff#1\lopoff#1#2}
\long\def\lopoff\\#1#2\lopoff#3#4{\def#4{#1}\def#3{#2}}

\def\ismember#1\of#2{\foundfalse{\let\given=#1%
    \def\\##1{\def\next{##1}%
    \ifx\next\given{\global\foundtrue}\fi}#2}}

% Les commandes
\def\section#1{\vskip1truecm
               \global\def\currenvir{section}
               \global\advance\secno by1\global\prmno=0
               {\bf \number\secno. {#1}}
              \vglue3pt}

\def\subsection{\global\def\currenvir{subsection}
                \global\advance\prmno by1
                \ind{ (\number\secno.\number\prmno) }}
\def\subsec{\global\def\currenvir{subsection}
                \global\advance\prmno by1
                { (\number\secno.\number\prmno)\ }}

\def\proclaim#1{\global\advance\prmno by 1
                {\bf #1 \the\secno.\the\prmno$.-$ }}

\long\def\th#1 \enonce#2\endth{%
   \medbreak\proclaim{#1}{\it #2}\global\def\currenvir{th}
\smallskip}

\def\bib#1{\rm #1}
\long\def\thr#1\bib#2\enonce#3\endth{%
\medbreak{\global\advance\prmno by 1\bf#1\the\secno.\the\prmno\ 
\bib{#2}$\!.-$ } {\it
#3}\global\def\currenvir{th}\smallskip}
%usage:\thr Proposition \bib Dupont \enonce

\def\rem#1{\global\advance\prmno by 1
{\it #1} \the\secno.\the\prmno$.-$ }

% CROSS-REFERENCES
\def\isinlabellist#1\of#2{\notfoundtrue%
   {\def\given{#1}%
    \def\\##1{\def\next{##1}%
    \lop\next\to\za\lop\next\to\zb%
    \ifx\za\given{\zb\global\notfoundfalse}\fi}#2}%
    \ifnotfound{\immediate\write16%
                 {Warning - [Page \the\pageno] 
{#1} No reference found}}%
 \fi}%
\def\ref#1{\ifx\labellist\empty{\immediate\write16
                 {Warning - No references found at all.}}
               \else{\isinlabellist{#1}\of\labellist}\fi}

\def\newlabel#1#2{\rightappenditem{\\{#1}\\{#2}}\to\labellist}
\def\labellist{}
\def\label#1{%
  \def\given{th}%
  \ifx\given\currenvir%
    {\hwrite\lbl{\string\newlabel{#1}{\number\secno.\number\prmno}}}\fi%
  \def\given{section}%
  \ifx\given\currenvir%
    {\hwrite\lbl{\string\newlabel{#1}{\number\secno}}}\fi%
  \def\given{subsection}%
  \ifx\given\currenvir%
    {\hwrite\lbl{\string\newlabel{#1}{\number\secno.\number\prmno}}}\fi%
  \def\given{subsubsection}%
  \ifx\given\currenvir%
  {\hwrite\lbl{\string%
    \newlabel{#1}{\number\secno.\number\subsecno.\number\subsubsecno}}}\fi
  \ignorespaces}
\newwrite\lbl
\def\begin{\newlabel{lambda}{1.1}
\newlabel{chern}{1.2}
\newlabel{psi}{1.4}
\newlabel{split}{1.5}
\newlabel{inv}{1.6}
\newlabel{R(T)}{2.1}}
\magnification 1250
\pretolerance=500 \tolerance=1000  \brokenpenalty=5000
\mathcode`A="7041 \mathcode`B="7042 \mathcode`C="7043
\mathcode`D="7044 \mathcode`E="7045 \mathcode`F="7046
\mathcode`G="7047 \mathcode`H="7048 \mathcode`I="7049
\mathcode`J="704A \mathcode`K="704B \mathcode`L="704C
\mathcode`M="704D \mathcode`N="704E \mathcode`O="704F
\mathcode`P="7050 \mathcode`Q="7051 \mathcode`R="7052
\mathcode`S="7053 \mathcode`T="7054 \mathcode`U="7055
\mathcode`V="7056 \mathcode`W="7057 \mathcode`X="7058
\mathcode`Y="7059 \mathcode`Z="705A
\def\spacedmath#1{\def\packedmath##1${\bgroup\mathsurround =0pt##1\egroup$}
\mathsurround#1
\everymath={\packedmath}\everydisplay={\mathsurround=0pt}}
\def\nospacedmath{\mathsurround=0pt
\everymath={}\everydisplay={} } \spacedmath{2pt}

\def\phfl#1#2{\normalbaselines{\baselineskip=0pt
\lineskip=10truept\lineskiplimit=1truept}\nospacedmath\smash {\mathop{\hbox to
8truemm{\rightarrowfill}}
\limits^{\scriptstyle#1}_{\scriptstyle#2}}}
\def\hfl#1#2{\normalbaselines{\baselineskip=0truept
\lineskip=10truept\lineskiplimit=1truept}\nospacedmath\smash{\mathop{\hbox to
12truemm{\rightarrowfill}}\limits^{\scriptstyle#1}_{\scriptstyle#2}}}
\def\diagram#1{\def\normalbaselines{\baselineskip=0truept
\lineskip=10truept\lineskiplimit=1truept}   \matrix{#1}}
\def\wfl#1#2{\llap{$\scriptstyle#1$}\left\uparrow\vbox
to 6truemm{}\right.\rlap{$\scriptstyle#2$}}
\def\vfl#1#2{\llap{$\scriptstyle#1$}\left\downarrow\vbox
to 6truemm{}\right.\rlap{$\scriptstyle#2$}}

\def\mono{\lhook\joinrel\mathrel{\longrightarrow}}
\def\iso{\vbox{\hbox to .8cm{\hfill{$\scriptstyle\sim$}\hfill}
\nointerlineskip\hbox to .8cm{{\hfill$\longrightarrow $\hfill}} }}

\def\sdir_#1^#2{\mathrel{\mathop{\kern0pt\oplus}\limits_{#1}^{#2}}}
\def\pprod_#1^#2{\raise
2pt \hbox{$\mathrel{\scriptstyle\mathop{\kern0pt\prod}\limits_{#1}^{#2}}$}}

\font\eightrm=cmr8         \font\eighti=cmmi8
\font\eightsy=cmsy8        \font\eightbf=cmbx8
\font\eighttt=cmtt8        \font\eightit=cmti8
\font\eightsl=cmsl8        \font\sixrm=cmr6
\font\sixi=cmmi6           \font\sixsy=cmsy6
\font\sixbf=cmbx6\catcode`\@=11
\def\eightpoint{%
  \textfont0=\eightrm \scriptfont0=\sixrm \scriptscriptfont0=\fiverm
  \def\rm{\fam\z@\eightrm}%
  \textfont1=\eighti  \scriptfont1=\sixi  \scriptscriptfont1=\fivei
  \def\oldstyle{\fam\@ne\eighti}\let\old=\oldstyle
  \textfont2=\eightsy \scriptfont2=\sixsy \scriptscriptfont2=\fivesy
  \textfont\itfam=\eightit
  \def\it{\fam\itfam\eightit}%
  \textfont\slfam=\eightsl
  \def\sl{\fam\slfam\eightsl}%
  \textfont\bffam=\eightbf \scriptfont\bffam=\sixbf
  \scriptscriptfont\bffam=\fivebf
  \def\bf{\fam\bffam\eightbf}%
  \textfont\ttfam=\eighttt
  \def\tt{\fam\ttfam\eighttt}%
  \abovedisplayskip=9pt plus 3pt minus 9pt
  \belowdisplayskip=\abovedisplayskip
  \abovedisplayshortskip=0pt plus 3pt
  \belowdisplayshortskip=3pt plus 3pt 
  \smallskipamount=2pt plus 1pt minus 1pt
  \medskipamount=4pt plus 2pt minus 1pt
  \bigskipamount=9pt plus 3pt minus 3pt
  \normalbaselineskip=9pt
  \setbox\strutbox=\hbox{\vrule height7pt depth2pt width0pt}%
  \normalbaselines\rm}\catcode`\@=12
\font\tensan=cmssdc10
\font\sevensan=cmssdc10 at 7pt
\font\fivesan=cmssdc10 at 5pt
\newfam\sanfam
\textfont\sanfam=\tensan \scriptfont\sanfam=\sevensan
  \scriptscriptfont\sanfam=\fivesan
  \def\san{\fam\sanfam\tensan}%

\newcount\noteno
\noteno=0
\def\up#1{\raise 1ex\hbox{\sevenrm#1}}
\def\note#1{\global\advance\noteno by1
\footnote{\parindent0.4cm\up{\number\noteno}\
}{\vtop{\eightpoint\baselineskip12pt\hsize15.5truecm\noindent
#1}}\parindent 0cm}
%\font\san=cmssdc10
\def\ext{\hbox{\san \char3}}
\def\sym{{\san \char83}}
\def\car{{\san \char88}}
\font\gragrec=cmmib10
\def\lam{\hbox{\gragrec \char21}}

\def\pc#1{\tenrm#1\sevenrm}
\def\tx{\kern-1.5pt -}
\def\cqfd{\kern 2truemm\unskip\penalty 10000\vrule height
4pt depth 0pt width 4pt\medbreak} 
\def\no{n\up{o}\kern 2pt}
\def\ind{\par\hskip 0.8truecm\relax}
\def\indp{\par\hskip 0.4truecm\relax}

\def\rond{\kern 1pt{\scriptstyle\circ}\kern 1pt}
\def\Ker{\mathop{\rm Ker}\nolimits}

\def\det{\mathop{\rm det}\nolimits}

\def\dim{\mathop{\rm dim}\nolimits}

\def\Tr{\mathop{\rm Tr}\nolimits}
\frenchspacing
\vsize = 25truecm
\hsize = 16truecm
%\hoffset = -.15truecm
\voffset = -.5truecm
\parindent=0cm
\baselineskip15pt
\def\pr{{\it Proof} : }
\begin\null\vskip1mm
\centerline{\bf  Chern classes for representations
 of  reductive groups}
\smallskip
\smallskip \centerline{Arnaud {\pc BEAUVILLE}} 
\vskip1.2cm 
{\bf Introduction}
\smallskip
\ind
Let $G$ be a complex connected reductive group, and let
$R(G)$ be its representation ring. As an abelian group
$R(G)$ is spanned by the finite-dimensional
representations $[V]$ of $G$, with the relations
$[V\oplus W]=[V]+[W]$; the ring structure is defined by
the tensor product. Moreover the exterior product  of
representations give rise to a sequence of operations
$\lambda^p:R(G)\rightarrow R(G)$ which make
$R(G)$ into a {\it $\lambda$\tx ring} -- see
(\ref{lambda}) below for the definition. 
\ind In the course on his work on the
Riemann-Roch theorem, Grothendieck had the remarkable
insight that this purely algebraic structure is enough
to define Chern classes, without any reference to a
cohomology or Chow ring. He associated to any
$\lambda$\tx ring $R$ a filtration of $R$, the
$\gamma$\tx {\it filtration}$\!$;  the Chern classes take
value in the associated graded ring
${\rm gr}\,R$, and the Chern character  is a ring
homomorphism
${\rm ch}:R\rightarrow \pprod_p^{}\ {\rm gr}^p_{\bf Q}$,
where ${\rm gr}_{\bf Q}\,R={\rm gr}\,R\otimes{\bf Q}$.
When applied to the Grothendieck ring $K(X)$ of vector
bundles on a smooth algebraic variety $X$ these
definitions give back the classical ones, at least
modulo torsion: the graded ring
${\rm gr}_{\bf Q}\,K(X)$ coincides with the Chow ring
$CH(X)_{\bf Q}$, and the Chern classes with the usual
ones.
\ind The aim of this note is to compute the Chern
classes for the representation ring $R(G)$ -- a simple
exercise which I have been unable to find in the
literature. Let
${\goth g}$ be the Lie algebra of $G$; we will assume
that $G$, and therefore ${\goth g}$, are defined over
${\bf Q}$ (alternatively, we could  without any loss take
our Chern classes in ${\rm gr}_{\bf C}\,R(G)$). We
denote by
${\rm Pol}({\goth g})^{\rm inv}$ the ring of polynomial
functions on ${\goth g}$ which are 
invariant under the action of the adjoint group. 
\medbreak
{\bf Theorem}$.-$ a) {\it The graded ring ${\rm
gr}_{\bf Q}\,R(G)$ is canonically isomorphic to the
ring ${\rm Pol}({\goth g})^{\rm inv}$ of invariant
polynomial functions on ${\goth g}$.}
\ind b) {\it Let $\rho$ be a representation of $G$, and
$L\rho$ the corresponding representation of ${\goth g}$.
The total Chern class $c(\rho)$ is equal to the
invariant function $\det(1+L\rho)$, and the Chern
character ${\rm ch}(\rho)$ to} $\Tr(\exp L\rho)$.
\smallskip 
\ind We have of course an explicit description of the
ring $R(G)$. Let $T$ be a maximal torus of $G$; the Weyl
group $W$ acts on the ring $R(T)$, and the restriction
map $R(G)\rightarrow R(T)$ identifies $R(G)$ with the
invariant sub-ring $R(T)^W$. As an abelian group $R(T)$
is spanned by one-dimensional elements (we will say that
it is a {\it split} $\lambda$\tx ring), which makes easy
to compute its $\gamma$\tx filtration and Chern
classes. The theorem follows easily once we know that
the $\gamma$\tx filtration of $R(T)$ induces that of
$R(G)$. This  turns out to be  a general fact for
invariant sub-rings of split
$\lambda$\tx rings; we will deduce it from the behaviour
of the $\gamma$\tx filtration under the Adams operations
(Proposition \ref{inv}).
\ind In section 3 we discuss an application. Let $P$ be a
principal
$G$\tx bundle over a base $B$ (which may be a variety,
or an arbitrary topos); the Theorem provides a simple
definition of the characteristic classes of $P$ in the
graded ring
${\rm gr}_{\bf Q}\,K(B)$, and a simple way of computing 
the Chern classes of the associated vector bundles.

\section{Generalities on $\lam$\tx rings}
\ind In this section we recall the definition and basic
properties of $\lambda$\tx rings. Standard references
are [SGA6] or [F-L]; we follow the terminology of [SGA6],
Expos\'e V.\smallskip 
\subsection\label{lambda} A $\lambda$\tx {\it ring} $R$
is a commutative ring with two more pieces of structure:
\indp -- An augmentation, that is a 
ring homomorphism $\varepsilon:R\rightarrow {\bf Z}$;
\indp -- A $\lambda$\tx structure, that is a sequence of maps
$\lambda^i:R\rightarrow R$ such that, for any $x,y$ in $R$
and $n\in{\bf N}$,
$$\lambda^0(x)=1\qquad \lambda^1(x)=x \qquad
\lambda^n(x+y)=\sum_{p+q=n}\lambda^p(x)\,\lambda^q(y)\ .$$
If we put
$\displaystyle
\lambda_t(x)=\sum_p\lambda^p(x)\,t^p\in R[[t]]$, the
last condition is equivalent to
$$\lambda_t(x+y)=\lambda_t(x)\,\lambda_t(y)\ .$$

\ind Moreover we want formulas giving 
$\lambda^p(xy)$  and $\lambda^p(\lambda^q(x))$ as 
polynomials in $\lambda^1(x),\ldots,\lambda^p(x);
\lambda^1(y),\ldots,\lambda^p(y)$ and 
$\lambda^1(x),\ldots,\lambda^{pq}(x)$ respectively. A
convenient way of expressing these is to introduce a
sequence $(\psi^k)_{k\ge 1}$ of additive endomorphisms of $R$, 
the {\it Adams operations}, defined by 
$$\sum_{k\ge 1}\psi^k(x)\,(-t)^{k-1}=\lambda_{t}(x)^{-1} \ 
{d\over dt}\lambda_{t}(x)\ .$$
\ind Then the condition on the $\lambda$\tx structure
means that the
$\psi^k$ are ring endomorphisms, and satisfy $\psi^k\rond
\psi^\ell=\psi^{k\ell }$ for $k,\ell \ge 1$.

\subsection\label{chern}   Grothendieck associates to
this situation a second $\lambda$\tx structure, defined
by
$\gamma_t(x)=\lambda_{t\over 1-t}(x)$, and a
decreasing filtration of $R$, the $\gamma$\tx {\it
fitration} $(\Gamma ^p(R))_{p\ge 0}$: we put
$\Gamma ^0=R$, $\Gamma ^1=\Ker \varepsilon$, and $\Gamma ^p$
is spanned by the elements
$\gamma^{i_1}(x_1)\ldots\gamma^{i_k}(x_k)$ with
$x_1,\ldots,x_k\in\Gamma ^1$ and $i_1+\cdots+i_k=p$. Let
${\rm gr}\, R=\sdir_{p\ge 0}^{}\ \Gamma ^p/\Gamma^{p+1}$
be the associated graded ring. The {\it Chern classes}
$(c_p(x))_{p\ge 0}$ of an element
$x\in R$ are defined by
$$c_p(x)=\gamma^p(x-\varepsilon(x))\quad {\rm in}\quad
\Gamma ^p/\Gamma^{p+1}\ .$$
 
\subsection   A crucial point in what follows will be
the behaviour of the $\gamma$\tx filtration with
respect to the Adams operations.
Let us say that an element $x$ of $R$ has 
$\lambda$\tx {\it dimension} $n$ if it satisfies
$\varepsilon(x)=n$ and $\lambda^i(x)=0$ for $i>n$; we say
that $R$ is {\bf Q}-{\it finite
$\!\lambda$\tx dimensional}  if the ${\bf Q}$\tx vector
space
$R\otimes{\bf Q}$ is spanned by  elements of finite
$\lambda$\tx dimension. For such a
$\lambda$\tx ring we have, for each  $k\ge 1$ and $x\in\Gamma
^p(R)$: 
$$\psi^k(x)\equiv k^px\quad\hbox{mod. }{\bf Q}\Gamma
^{p+1}(R)\leqno{\subsec}$$
([F-L], III, Proposition 3.1).\label{psi}
\smallskip \subsection\label{split} Let us say that the
$\lambda$\tx ring
$R$ is {\it split} if is generated as an abelian group
by elements of $\lambda$\tx dimension 1. In this case we
have $\Gamma^p(R)={\goth r}^p $, where
${\goth r}$ is the augmentation ideal of $R$: indeed we
have ${\goth r}^p\i\Gamma ^p$ because $\gamma^1$ is the
identity, and $\gamma_t(x)\in\sum_{p\ge 0}{\goth
r}^p\,t^p$ for all $x\in {\goth r}$ because this holds
for $x+y$ if it holds for $x$ and $y$, and
$\gamma_t(\xi-1)=1+(\xi-1)t$ if $\xi$ has $\lambda$\tx
dimension 1.
\smallskip 
\th Proposition
\enonce Let $R$ be a split $\lambda$\tx ring, and $G$ a
finite group of automorphisms of the $\lambda$\tx ring
$R$. Then the $G$\tx invariant elements form a 
 sub-$\lambda$\tx ring $S$ of
$R$. The  filtrations  $(\Gamma ^p(S))_{p\ge 0}$ and $(\Gamma
^p(R)\cap S)_{p\ge 0}$ coincide in $S\otimes{\bf Q}$.
\endth\label{inv}
\pr Since the action of $G$ commutes with the
$\lambda^i$ the $\lambda$\tx structure of $R$ induces a 
$\lambda$\tx structure on $S$. Moreover the ${\bf Q}$\tx 
vector space $S\otimes{\bf Q}$ is spanned  by the elements
$\sum_{g\in G}g\xi$, where
$\xi$ is an element of $R$ of $\lambda$\tx dimension 1;
these elements have $\lambda$\tx dimension $|G|$, thus
$S$ is  ${\bf Q}$\tx finite $\lambda$\tx dimensional.

\ind To alleviate the notation we write $R,S$ instead of 
$R\otimes{\bf Q}$ and $S\otimes{\bf Q}$, and $\Gamma
^p(R),\Gamma ^p(S)$ instead of ${\bf Q}\Gamma ^p(R),{\bf
Q}\Gamma ^p(S)$. We have $\Gamma^p(S)\i\Gamma ^p(R)\cap S$
for each
$p$. Let us first
prove that the two filtrations define the same topology. 
Let ${\goth r}=\Gamma ^1(R)$ be the augmentation ideal of
$R$, and ${\goth s}=S\cap{\goth r}$. Since $R$ is a
finite $S$\tx module, the ring $R/{\goth s}R$ is artinian; 
thus there exists an integer
$\nu$ such that ${\goth r}^\nu\i {\goth s}R$, and therefore 
${\goth r}^{n\nu}\i {\goth s}^nR$ for all $n$. 

\ind As $S$\tx modules $S$ is a direct summand of  $R$
(consider the projector $r\mapsto$ $\displaystyle
{1\over |G|}\sum_{g\in G}gr$). Thus for every ideal
${\goth a}$ of $S$, the induced homomorphism $S/{\goth
a}\rightarrow R/{\goth a}R$ is injective, which means
${\goth a}R\cap S={\goth a}$. Using (\ref{split}) we
obtain
$$\Gamma ^{p\nu}(R)\cap S={\goth r}^{p\nu}\cap S\i
{\goth s}^pR\cap S
= {\goth s}^{p}\i \Gamma^{p}(S)\ .$$
\ind Let us prove now the
 inclusion $\Gamma ^p(R)\cap S\i \Gamma^p(S)$ by
induction on $p$, the case
$p=0$ being obvious. We have just seen
that there exists an integer
$N$ such that $\Gamma ^N(R)\cap S\i \Gamma ^p(S)$; let
$n$ be the smallest integer with that property. Assume
$n>p$. Let $x\in
\Gamma ^{n-1}(R)\cap S$, and let $k\ge 2$ be an integer;
by (\ref{psi}) we have
$$\psi^k(x)\equiv k^{n-1}x\quad {\rm (mod.}\
\Gamma ^{n}(R)\cap S)\ .$$ \ind On the other hand we
have
$x\in \Gamma ^{p-1}(S)$ by the induction hypothesis;
since $S$ is ${\bf Q}$\tx finite
$\lambda$\tx dimensional (\ref{psi}) gives
$$\psi^k(x)\equiv k^{p-1}x \quad {\rm (mod.}\
\Gamma^p(S))\ .$$ 
\ind Since $\Gamma ^n(R)\cap S\i \Gamma ^p(S)$ we get
$x\in \Gamma ^p(S)$ for all $x\in \Gamma ^{n-1}(R)\cap
S$, contra\-dicting the choice of $n$. Therefore we have
$n\le p$, hence
$\Gamma ^p(R)\cap S=\Gamma^p(S)$.\cqfd

\section{The $\lam$\tx ring R(G)}
\subsection \label{R(T)}Let $T$ be a connected
multiplicative group, and $\car$ its character group;
this is a free finitely generated abelian group. The
representation ring  $R(T)$ is isomorphic to the group
algebra ${\bf Z}[\car]$; we denote by
$([\alpha])_{\alpha\in\car}$ its canonical basis. The
augmentation and the $\lambda$\tx structure are
characterized by $\varepsilon([\alpha])=1$ and
$\lambda_t([\alpha])=1+t\alpha$. 
\ind Let
${\goth r}$ be the augmentation ideal in $R(T)$. We
have a canonical isomorphism
$\car\rightarrow {\goth r}/{\goth r}^2$ which maps a
character $\alpha$ to the class of $[\alpha]-1$  (observe
that
$[\alpha\beta]-1\equiv$ $([\alpha]-1)+([\beta]-1)$
(mod.$\,{\goth r}^2$). Since the rings $R(T)$ and
$R(T)/{\goth r}\cong {\bf Z}$ are regular, the canonical
map
$\sym({\goth r}/{\goth r}^2)\rightarrow$
$\sdir_{p\ge 0}^{}{\goth r}^p/{\goth r}^{p+1}$ is
bijective; using (\ref{split}) we get
a canonical isomorphism
$$\varphi : \sym(\car)\iso {\rm gr}\, R(T)\ .$$
\ind Under this isomorphism the element
$\gamma^1([\alpha]-1)$ of ${\goth r}/{\goth r}^2$
corresponds to $\alpha\in \car\i\sym(\car)$. Thus the
total Chern class $c(\rho)$ of an element $\rho=
[\alpha_1]+\cdots+[\alpha_d]$ is given by
$c(\rho)=\pprod_i^{}(1+\alpha_i)$. 

 \ind  Let ${\goth h}$ be the Lie algebra of $T$,
viewed as a vector space over ${\bf Q}$. We have a
canonical isomorphism $\car\otimes{\bf Q}\rightarrow
{\goth h}^*$, which associates to a character
$\alpha:T\rightarrow {\bf G}_m$ its derivative $L\alpha:{\goth
h}\rightarrow {\bf Q}$. Thus we can identify 
$\sym(\car)\otimes{\bf Q}$ to the algebra ${\rm Pol}({\goth
h})$ of polynomial maps ${\goth h}\rightarrow {\bf
Q}$; the above formula becomes
$c(\rho)=\pprod_i^{}(1+L\alpha_i)$ in ${\rm Pol}({\goth
h})$.

\ind For $\alpha\in\car$, let  ${\bf
C}_\alpha$ denote the one-dimensional representation with
character $\alpha$. The element
$\rho=[\alpha_1]+\ldots+[\alpha_d]$ of $R(T)$ is the
class of the representation $V={\bf
C}_{\alpha_1}\oplus\ldots\oplus {\bf C}_{\alpha_d}$. In the
corresponding representation $L\rho$ of ${\goth h}$, an
element $H$ of ${\goth h}$ acts through the diagonal
matrix ${\rm diag}\,(L\alpha_1(H),\ldots,L\alpha_d(H))$.
Thus $c(\rho)$ is equal to the function $\det(1+L\rho)$ in 
${\rm Pol}({\goth h})$.
\ind The Chern character gives an homomorphism ${\rm
ch}$ of $R(T)$ into the ring $\widehat{\rm Pol}({\goth h})$
of formal series on ${\goth h}$. We have ${\rm
ch}([\alpha])=e^{c_1([\alpha])}=e^{ L\alpha}$, hence
$${\rm ch}(\rho)=e^{L\alpha_1}+\ldots+e^{
L\alpha_d}=\Tr(\exp L\rho)\quad\hbox{in }\widehat{\rm
Pol}({\goth h})\ .$$ \smallskip 
\rem{Remark} The exponential  morphism of
formal groups $\exp:\widehat{\goth h}\rightarrow
\widehat{T}$ induces an injective homomorphism
$\exp^*:{\cal O}_{\widehat{T},1}={\bf
Q}[[\car]]\rightarrow
{\cal O}_{\widehat{\goth h},0}=\widehat{\rm Pol}({\goth
h})$, which maps $[\alpha]$ to $e^{L\alpha}$. The Chern
character is the composition of this map with the
injection $R(T)={\bf Z}[\car]\mono{\bf Q}[[\car]]$.
\medskip 
\subsection  Let $G$ be a complex connected reductive
group, ${\goth g}$ its Lie algebra,
$T$ a maximal torus of $G$ and ${\goth h}$ its Lie
algebra; we can assume that $G$,
$T$, ${\goth g}$ and ${\goth h}$ are  defined over ${\bf
Q}$. The Weyl group $W$ acts on
$T$, hence on the ring $R(T)$; restriction to
$T$ induces a homomorphism of $\lambda$\tx rings
$R(G)\rightarrow R(T)$, whose image is the invariant
sub-ring $R(T)^W$ ([SGA6], Expos\'e 0 App., Th. 1.1). By
Proposition
\ref{inv}, the graded ring ${\rm gr}_{\bf Q}\,R(G)$ is
isomorphic to 
$({\rm gr}_{\bf Q}\,R(T))^W$, that is to the ring 
${\rm Pol}({\goth h})^W$ of invariant polynomials on
${\goth h}$ (\ref{R(T)}).
\ind Let $\rho$ be a representation  of $G$. From the
commutative diagrams
$$\diagram{R(T)&\hfl{c}{}&\kern-2mm\widehat{\rm Pol}({\goth
h})\cr
\wfl{}{}&&\kern-2mm\wfl{}{}&&\cr
R(G)&\hfl{c}{}&\widehat{\rm Pol}({\goth h})^W}\qquad\qquad
\diagram{R(T)&\hfl{\rm ch}{}&\kern-2mm\widehat{\rm
Pol}({\goth h})\cr
\wfl{}{}&&\kern-2mm\wfl{}{}\cr
R(G)&\hfl{\rm ch}{}&\widehat{\rm Pol}({\goth
h})^W}$$
we obtain
$$c(\rho)=\det(1+L\rho)\ \ \hbox{in }\ {\rm Pol}({\goth
h})^W\ ;\quad{\rm ch(\rho)=\Tr(\exp L\rho)}\ \ \hbox{in }\ 
\widehat{\rm Pol}({\goth h})^W\ .$$
 The restriction map
${\rm Pol}({\goth g})\rightarrow {\rm Pol}({\goth h})$
induces an isomorphism ${\rm Pol}({\goth
g})^{inv}\rightarrow {\rm Pol}({\goth h})^W$ ([B], \S 8,
Th. 1). Since the functions $\det(1+L\rho)$ and
$\Tr(\exp L\rho)$ are invariant, the above equalities
hold as well in ${\rm Pol}({\goth g})^{inv}$ and
$\widehat{\rm Pol}({\goth g})^{inv}$ respectively. This
proves the theorem stated in the introduction.

\section{Application: Chern classes of associated
bundles}
\subsection Let $B$ be an algebraic variety, and $P$ a
principal $G$\tx bundle over $B$ (for what follows $B$
could be as well a scheme, or even a topos). To any
representation $\rho:G\rightarrow  GL(V)$ is associated
 a vector bundle
$P^{\rho}=P\times^GV$ on $B$; we define in this way a
homomorphism of
$\lambda$\tx rings $b_P^*:R(G)\rightarrow K(B)$\note{In
fancy terms, the principal bundle $P$ corresponds to a
morphism $b_P$ of $B$ into the classifying topos
(or stack) $BG$, and $b_P^*$ is just the pull-back map.}.
In view of the isomorphism described above, it induces a
homomorphism of graded rings
$$c^{}_P:{\rm Pol}({\goth g})^{inv}\rightarrow {\rm
gr}_{\bf Q}K(B)\ ,$$called the characteristic
homomorphism. 
\ind Let $\ell =\dim T$. Recall that there are
homogeneous functions $I_1,\ldots,I_\ell$ in ${\rm
Pol}({\goth g})^{inv}$ such that ${\rm
Pol}({\goth g})^{inv}={\bf Q}[I_1,\ldots,I_\ell]$ ([B],
\S 8, Th\'eor\`eme 1). The elements $c_P^{(i)}:=c^{}_P(I_i)$
may be called the {\it characteristic classes} of $P$
(but note that they depend on the choice of the
generating sequence
$(I_1,\ldots,I_\ell)$). If $B$ is a smooth variety, the
graded ring ${\rm gr}_{\bf Q}\,K(B)$ is canonically
isomorphic to the rational Chow ring $CH_{\bf Q}(B)$
([SGA6], Expos\'e XIV, \no 4); our characteristic classes
correspond under this isomorphism to those defined in
[V] and [E-G].
\th Proposition
\enonce Let $\rho$ be a representation of $G$; write
$\det(1+L\rho)=$ $F(I_1,\ldots,I_\ell )$, where $F$ is a
polynomial in $\ell $ indeterminates. Let $P$ be a
principal $G$\tx bundle on $B$, with characteristic
classes $c_P^{(1)},\ldots,c_P^{(\ell)}$ in ${\rm
gr}_{\bf Q}\,K(B)$. Then the total Chern class in ${\rm
gr}_{\bf Q}\,K(B)$ of the
associated bundle $P^\rho$  is
$$c(P^\rho)=F(c_P^{(1)},\ldots,c_P^{(\ell)})\ .$$
\ind Similarly, if $\Tr(\exp
L\rho)=G(I_1,\ldots,I_\ell)$, with
$G\in {\bf Q}[[T_1,\ldots,T_\ell]]$, we have \break
${\rm ch}(P^\rho)=G(c_P^{(1)},\ldots,c_P^{(\ell)})$.
\endth\smallskip 
\pr This follows from the Theorem, the  commutative
diagram
\vskip-.5cm$$\diagram{R(G)&\hfl{b_P^*}{}&K(B)\cr
\vfl{c}{} &&\vfl{}{c}\cr \widehat{\rm Pol}({\goth g})^{inv}&
\hfl{c_P^{}}{}&\widehat{\rm gr}_{\bf Q}\,K(B)}$$\vskip-.2cm
and the corresponding diagram for the Chern
character.\cqfd
\rem{Examples} The Proposition provides a way of
computing the characteristic classes in terms of Chern
classes, at least for  classical groups. We use the
standard generating system
$(I_1,\ldots,I_\ell )$ given for instance in [B], \S 13.
We denote by $E_P$ the vector bundle on $B$ associated to
$P$ through the standard representation of
$GL(n),\,SO(n)\ \hbox{or }{\rm Sp}(n)$ in ${\bf C}^n$.
\ind  a) For
$G=GL(\ell)$, we  have $I_p(A)=\Tr (\ext^pA)$ for
$A\in{\goth g}{\goth l}(\ell )$ and $1\le p\le \ell $;
this gives $c_P^{(p)}=c_p(E_P)$ for  $1\le p\le
\ell $.
\ind b) For $G={\rm Sp}(2\ell )$ or $SO(2\ell +1)$, we
have $I_p(A)=\Tr (\ext^{2p}A)$; this gives 
$c_P^{(p)}=c_{2p}(E_P)$ for $1\le p\le \ell $ (the
vector bundle $E_P$ is isomorphic to its dual, thus its
odd Chern classes with rational coefficients vanish).
\ind c) For $G=SO(2\ell)$, we realize ${\goth g}$ as the
space of skew-symmetric matrices in ${\goth g}{\goth
l}(2\ell )$; we have $I_p(A)=\Tr (\ext^{2p}A)$ for
$p<\ell
$, and
$I_{2\ell }(A)={\rm Pf}(A)$. We obtain
$c_P^{(p)}=c_{2p}(E_P)$ for $1\le p\le \ell -1$, and
$c_P^{(\ell )}$ is a class in ${\rm gr}^{\ell }_{\bf
Q}\,K(B)$ with square $c_{2\ell }(E_P)$.

\vskip2cm\frenchspacing
\centerline{ REFERENCES} \vglue15pt\baselineskip12.8pt
\def\num#1{\smallskip\item{\hbox to\parindent {\enskip
[#1]\hfill}}}\parindent=1.3cm 
\num{B} N. {\pc BOURBAKI}: {\sl Groupes et alg\`ebres de
Lie}, chapitre VIII. Hermann, Paris (1975).
\num{E-G} D. {\pc EDIDIN}, W. {\pc GRAHAM}: {\sl 
Characteristic classes in the Chow ring}. J. Algebraic
Geom. {\bf 6} (1997), 431--443. 
\num{F-L} W. {\pc FULTON}, S.
{\pc LANG}: {\sl Riemann-Roch algebra}.  Grundlehren der
Math.
 {\bf 277}.  Springer-Verlag, New York (1985).
\num{SGA6} {\sl Th\'eorie des intersections et th\'eor\`eme de
Riemann-Roch}.  S\'eminaire de G\'eom\'etrie Alg\'ebrique
du Bois-Marie 1966--1967 (SGA 6). Dirig\'e par
P.~Berthelot, A. Grothendieck et L. Illusie.  Lecture
Notes in Math. {\bf  225}, Springer-Verlag, Berlin-New
York (1971). 
\num{V} A. {\pc VISTOLI}: 
  {\sl  Characteristic classes of principal bundles in algebraic intersection
theory}. Duke Math. J. {\bf 58} (1989),  299--315. 
\vskip1cm
\def\pc#1{\eightrm#1\sixrm}
\hfill\vtop{\eightrm\hbox to 5cm{\hfill Arnaud {\pc BEAUVILLE}\hfill}
 \hbox to 5cm{\hfill Laboratoire J.-A. Dieudonn\'e\hfill}
 \hbox to 5cm{\sixrm\hfill UMR 6621 du CNRS\hfill}
\hbox to 5cm{\hfill {\pc UNIVERSIT\'E DE} 
{\pc NICE}\hfill}
\hbox to 5cm{\hfill  Parc Valrose\hfill}
\hbox to 5cm{\hfill F-06108 {\pc NICE} Cedex 2\hfill}}
\end